\documentclass[preprint,11pt,authoryear,a4paper]{elsarticle}

\usepackage{adjustbox} 
\usepackage{setspace}
\usepackage{amssymb}
\usepackage{amsmath}
\usepackage{natbib}
\usepackage[ruled,vlined]{algorithm2e}
\usepackage{float}
\usepackage{graphicx}
\usepackage{subcaption}  
\usepackage{tabularx}    
\usepackage{siunitx}
\usepackage{threeparttable}
\usepackage{subcaption}  
\usepackage{makecell}    
\usepackage{booktabs}
\usepackage{amsthm}
\usepackage{amsmath}
\newtheorem{thm}{Theorem}

\newtheorem{corollary}[thm]{Corollary}

\newdefinition{rmk}{Remark}
\newproof{pf}{Proof}
\newproof{pot}{Proof of Theorem \ref{thm2}}

\usepackage[table]{xcolor}

\newcommand{\tablescale}{0.90} 

\newcommand{\rowht}{2.6ex} 
\newcommand{\rowstrut}{\rule{0pt}{\rowht}\rule[-0.8ex]{0pt}{0pt}}
\journal{European Journal of Operational Research}

\begin{document}

\begin{frontmatter}



\title{Sparse Convex Quantile Regression: A Generalized Benders Decomposition Approach}



\author[a]{Xiaoyu Luo\fnref{fn1}}
\author[a,b]{Chuanhou Gao\corref{cor1}%
}

\cortext[cor1]{Corresponding author. \newline  \hspace*{6mm}Email address: gaochou@zju.edu.cn}
\fntext[fn1]{This is the first author. \newline  \hspace*{6mm}Email address: 12135040@zju.edu.cn}

\affiliation[a]{organization={School of Mathematical Sciences, Zhejiang University},
            city={Hangzhou},
            postcode={310058}, 
            state={Zhejiang},
            country={China}}
            

\affiliation[b]{organization={Center for Interdisciplinary Applied Mathematics, Zhejiang University},
            city={Hangzhou},
            postcode={310058}, 
            state={Zhejiang},
            country={China}}
\begin{abstract}
We develop a scalable algorithmic framework for sparse convex quantile regression (SCQR), addressing key computational challenges in the literature. Enhancing the classical CQR model, we introduce $\ell_2$-norm regularization and an $\varepsilon$-insensitive zone to improve generalization and mitigate overfitting—both theoretically justified and empirically validated. Based on this extension, we improve the SCQR model and propose the first Generalized Benders Decomposition (GBD) algorithm tailored to this context, further strengthened by a novel local search-based Benders matheuristic. Extensive simulations and a real-world application to Sustainable Development Goals benchmarking demonstrate the accuracy, scalability, and practical value of our approach.
\end{abstract}


\begin{highlights}
\item We formulate convex quantile regression with \(\ell_2\)-regularization and adapt a primal cutting-plane method.
\item The proposed SCQR with \(\ell_2\)-regularization preserves the key \textit{quantile property}.
\item A generalized Benders decomposition algorithm is developed to solve the SCQR problem.
\item We develop a novel matheuristic that integrates local search into the Benders framework.
\end{highlights}


\end{frontmatter}


\section{Introduction}\label{sec:Intro}
Convex regression is a nonparametric technique used to estimate an unknown convex function from given data points. Unlike traditional linear regression, which assumes a linear relationship between input and output, convex regression relaxes this assumption and instead focuses on capturing the underlying convexity of the data. This approach is particularly useful in cases where the relationship between variables is inherently nonlinear but maintains an (approximate) convex structure, providing more flexibility while still preserving essential properties such as generalization and interpretability \citep{boyd2004convex}. Convex regression has gained significant attention due to its application in various fields, including economics, machine learning, and optimization, where capturing complex yet structured dependencies is crucial \citep{magnani2009convex, goldenshluger2006recovering, hannah2014semiconvex, topaloglu2003algorithm}. However, while convex regression offers significant flexibility, it is also prone to overfitting, particularly near the boundaries of training sample points, where the subgradients tend to grow excessively large \citep{liao2024convex}. This issue substantially undermines the generalization capacity of machine learning models. A common approach to alleviating this issue in the literature is to add a penalty to the objective loss function, such as the $\ell_2$-norm regularization \citep{bertsimas2021sparse, liao2024convex, mazumder2019computational}.\

Quantile regression is a statistical technique that extends classical linear regression by modeling the relationship between covariates and conditional quantiles of the response variable \citep{koenker1978regression}. Unlike ordinary least squares (OLS), which focuses on the conditional mean, quantile regression provides a fuller picture of the conditional distribution by estimating specific quantiles such as the median or other percentiles. This makes quantile regression particularly suitable for capturing heterogeneous effects, handling skewed distributions, and being robust to outliers \citep{koenker2001quantile}. The method minimizes the asymmetric quantile loss (pinball loss), which penalizes under- and over-estimations differently depending on the quantile level. As a result, it offers valuable insights into the impact of explanatory variables across the entire distribution of the outcome. Quantile regression is widely applied in various fields. For instance, financial analysts \citep{koenker2001quantile} may focus on extreme quantiles (e.g., 5th or 95th percentiles) to assess risk, while medical researchers may examine treatment effects across different risk groups \citep{yu2001bayesian}.
\

Recently, an increasing number of studies have explored Convex Quantile Regression (CQR) \citep{kuosmanen2015stochastic, wang2014nonparametric} and Convex Expectile Regression (CER) \citep{kuosmanen2021shadow, kuosmanen2020much}, which represent a promising integration of convex regression and quantile regression methodologies. By integrating these methodologies, CQR and CER allow for the estimation of conditional quantiles and expectiles while maintaining the convexity of the regression function. This combination enhances interpretability and ensures robustness in economic, financial, and operational research contexts involving nonlinear or asymmetric relationships\citep{kuosmanen2021shadow,dai2025optimal}. \cite{dai2023variable} introduced an $\ell_0$-constrained SCQR model and conducted a comparative study of its variable selection performance, benchmarked against $\ell_1$-norm regularization methods (Hastie, 2009). Through Monte Carlo simulations and an application to SDG performance evaluation across OECD countries, his results showed that the $\ell_0$-based approach better addresses the curse of dimensionality in high-dimensional settings.

However, limited research has addressed the development of scalable algorithms for solving the $\ell_0$-constrained SCQR problem. While \cite{dai2023variable} focused primarily on applying the SCQR framework in empirical analyses, including SDG benchmarking, the algorithmic aspects of solving such models efficiently remain underexplored. In this paper, we aim to bridge this gap by proposing the first decomposition-based algorithm tailored for the $\ell_0$-constrained SCQR problem. The main contributions of our work are outlined as follows:





\begin{itemize}
\item We address the CQR problem by incorporating an $\ell_2$-norm penalty on subgradients and the $\varepsilon$-insensitive zone, adapting the primal cutting-plane method from the literature \citep{bertsimas2021sparse, dai2023variable}, and demonstrate that the resulting SCQR model retains the fundamental quantile property.

\item We propose a GBD algorithm (Geoffrion, 1972) to solve the SCQR problem, representing the first scalable algorithm specifically designed for this purpose. Computational experiments show that the GBD algorithm delivers high-quality solutions within a few iterations. To further improve performance, we also develop a novel matheuristic that integrates local search with GBD.

\item Beyond computational aspects, we illustrate the practical value of SCQR through an application to the evaluation of SDG performance. By enabling frontier estimation at different quantile levels, our method captures heterogeneity in development performance and supports cross-country policy comparison, resource allocation, and strategic planning.
\end{itemize}

The structure of the paper is as follows. Section~\ref{sec:Prelimi} reviews the mathematical models for convex and quantile regression. Section~\ref{sec:CQR} introduces the convex quantile regression model with $\ell_2$-norm regularization and outlines the associated cutting-plane algorithm. Section~\ref{sec:Method} addresses the sparse convex quantile regression problem and presents the proposed generalized Benders decomposition method. Section~\ref{heuristic} develops a local search-based Benders matheuristic to improve incumbent solution quality. Finally, Section~\ref{sec:exp} reports computational results validating the effectiveness of the proposed approaches.

\subsection{Related literature}
There has been a growing body of research on decomposition algorithms and first-order optimization methods for variable selection, offering valuable insights into handling high-dimensional data and complex model structures. Since \cite{bertsimas2016best} introduced a mixed-integer optimization (MIO) framework with discrete first-order methods for best subset selection, exact sparse regression has seen renewed attention. \cite{bertsimas2020sparse} proposed a Benders-type dual cutting-plane method for sparse linear regression, and \cite{bertsimas2021sparse} extended this to sparse convex mean regression, leveraging smooth subproblems and well-structured duals for efficient solution.

Building on these foundations, \cite{chen2023sparse} addressed sparse linear quantile regression using MIO and first-order methods, while \cite{dai2023variable} introduced an $\ell_0$-constrained SCQR  solved via mixed-integer programming, showing the superiority of $\ell_0$ over $\ell_1$ regularization in high-dimensional settings. However, beyond the primal cutting-plane approach (CNLS-A) of \cite{dai2023variable}, scalable algorithmic frameworks for SCQR remain largely unexplored.

Motivated by these developments, we first enhance the classical CQR model by incorporating regularization techniques (Formulation (\ref{CQR})) inspired by support vector regression to mitigate overfitting. Building on this improved formulation, we propose the first generalized Benders decomposition algorithm for SCQR. In contrast to the sparse convex mean regression framework of \cite{bertsimas2021sparse}, where subproblems are smooth and symmetric, our  Benders subproblems (\ref{subproblem})  involve asymmetric, piecewise-linear (non-smooth) objectives, requiring rederivation of the dual and cut structures (see Theorem~\ref{min_to_max}). These structural differences lead to distinct algorithmic challenges in generating Benders cuts and ensuring convergence. To further enhance performance, we develop a novel improvement matheuristic that integrates
local search with the GBD algorithm, which is broadly applicable to general integer programming. From an application perspective, we demonstrate that our algorithm enables SCQR to successfully identify true variables at various quantile levels where sparse convex mean regression fails.


\section{Preliminaries}\label{sec:Prelimi}
In this section, we will formally introduce the mathematical formulation of convex regression and quantile regression.
\subsection{Convex regression}
Convex regression aims to estimate an unknown function \( f: \mathbb{R}^d \to \mathbb{R} \), where the observed response \( y \) can be expressed as: $y = f(\mathbf{x}) + \epsilon$, with the requirement that \( f \) is a convex function. Here, \( \mathbf{x} \in \mathbb{R}^d \) represents the predictor variables, and \(\epsilon\) is a random noise term that is assumed to have zero mean, i.e., \( \mathbb{E}[\epsilon] = 0 \). The convexity assumption of \( f \) implies that for any two points \( \mathbf{x}_1, \mathbf{x}_2 \in \mathbb{R}^d \) and any \( \lambda \in [0, 1] \), the following inequality holds:
$f(\lambda \mathbf{x}_1 + (1-\lambda) \mathbf{x}_2) \leq \lambda f(\mathbf{x}_1) + (1-\lambda)f(\mathbf{x}_2)$.


Given a set of observations \(\{(\mathbf{x}_i, y_i)\}_{i=1}^n\), where \( \mathbf{x}_i \in \mathbb{R}^d \) and \( y_i \in \mathbb{R} \), the goal of convex regression is to estimate \( f \) by minimizing the residual errors while ensuring the convexity of the estimated function. This formulation is infinite-dimensional, as the search space consists of continuous, real-valued convex functions. However, since the input data points are finite, the search space can be restricted to convex piecewise linear functions without any loss of accuracy\citep{boyd2004convex, kuosmanen2008representation}, thereby transforming the problem into a finite-dimensional one. The corresponding optimization problem \citep{boyd2004convex} can be written as: $\min_{\boldsymbol{\theta}, \boldsymbol{\beta}} \frac{1}{2} \sum_{i=1}^n \left( y_i - \theta_i \right)^2$, 
subject to the convexity constraints:$\theta_i + \boldsymbol{\beta}_i' (\mathbf{x}_j - \mathbf{x}_i) \leq \theta_j,  \forall i, j \in [n]$,
where \(\theta_i \) is the predicted response at \( \mathbf{x}_i \) and \(\boldsymbol{\beta}_i \in \mathbb{R}^d\) represents the subgradients of the estimated function $\hat{f}$ at \( \mathbf{x}_i \).

Given the optimal solutions \((\hat{\theta}, \hat{\boldsymbol{\beta}})\) to the above problem, we can reconstruct the explicit estimated function \(\hat{f}(\mathbf{x})\) as shown in \cite{kuosmanen2008representation}:
\begin{equation}
    \hat{f}(\mathbf{x}) = \max_{i=1, \dots, n} \left\{ \hat{\theta}_i + \hat{\boldsymbol{\beta}}_i' (\mathbf{x} - \mathbf{x_i}) \right\}.
\end{equation}

This function defines the estimated convex surface, which is a piecewise linear approximation determined by the observed data points. Each \(\hat{\boldsymbol{\beta}}_i\) represents a supporting hyperplane that characterizes the subgradient of the convex function \( f \) at the respective point \(\mathbf{x}_i \). In many applications, the true function $f$ may be concave; nonetheless, convex regression is still widely used. It is important to clarify that the regression function $f$ could either be globally convex or concave, depending on the sign of the convexity constraints (which can be reversed accordingly). In both cases, $f$ is the support of a convex set. Therefore, the term `convex regression' is used, as there are no concave sets in this context.

\subsection{Quantile Regression}
Quantile regression, introduced by \cite{koenker1978regression}, extends classical linear regression by estimating specific quantiles of the conditional distribution of the response variable \( y \) given the covariates \( \mathbf{x} \in \mathbb{R}^d \). Unlike ordinary least squares (OLS), which minimizes the squared loss to estimate the conditional mean, quantile regression estimates the conditional quantile \( \tau \in (0,1) \) of \( y \). Given a set of observations \(\{(\mathbf{x}_i, y_i)\}_{i=1}^n\), where \( \mathbf{x}_i \in \mathbb{R}^d \) and \( y_i \in \mathbb{R} \), the quantile regression model is formulated as:
\begin{equation}\label{model:QR}
    y_i = \mathbf{x}_i' \boldsymbol{\alpha} + \epsilon_i, \quad \text{with} \quad \mathbb{P}(\epsilon_i \leq 0 \mid \mathbf{x}_i) = \tau,
\end{equation}
where \(\boldsymbol{\alpha} \in \mathbb{R}^d \) is the vector of coefficients, and \( \tau \) represents the quantile level. The quantile regression estimator \(\hat{\boldsymbol{\alpha}}\) is obtained by solving the following optimization problem:
\begin{equation}\label{func:loss}
    \min_{\boldsymbol{\alpha} \in \mathbb{R}^d} \sum_{i=1}^n \rho_\tau \left( y_i - \mathbf{x}_i' \boldsymbol{\alpha} \right),
\end{equation}
In equation (\ref{model:QR}), the goal is to estimate the regression coefficients  \( \boldsymbol{\alpha} \) such that the conditional quantile \( \tau \) is correctly captured. Similar to how ordinary least squares (OLS) regression estimates the conditional mean by minimizing the squared loss function, quantile regression employs the \textit{check (pinball) loss function} \( \rho_{\tau} (y_i - \mathbf{x}_i' \boldsymbol{\alpha}) \), as formulated in equation (\ref{func:loss}). The check loss function is defined as:
\[
\rho_{\tau}(u) =
\begin{cases}
\tau u, & u \geq 0 \\
(\tau - 1)u, & u < 0
\end{cases}
\]

Then the conventional formulation of convex quantile regression can be formulated as follows \citep{dai2023variable}:
\begin{subequations}
    \label{conv_CQR}
    \begin{align}
         \min_{\boldsymbol{\beta}, \boldsymbol{\theta} , \xi, \xi^*} & \quad \sum_{i} \sum_{i=1}^{n} (\tau \xi_i + (1 - \tau) \xi_i^*) \\ \label{conv_positive}
    \text{s.t.} & \quad y_i - \theta_i \leq  \xi_i  \quad \forall i \in [n], \\ \label{conv_negative}
                & \quad \theta_i - y_i \leq \xi_i^* \quad \forall i \in [n], \\   \label{conv_convex cons}
                 &\quad \theta_i + \boldsymbol{\beta}_i'(\mathbf{x}_j-\mathbf{x}_i) \leq \theta_j \quad \forall i,j \in [n],\\
                & \quad \xi_i \geq 0, \quad \xi_i^* \geq 0 \quad \forall i \in [n].
    \end{align}
\end{subequations}
The model represents a conventional quantile regression formulation aiming to estimate the conditional quantile at a given level $\tau \in (0,1)$. In this formulation, $\theta_i$ denotes the estimated conditional quantile for sample $i$, while $\boldsymbol{\beta}_i$ represents the corresponding local linear coefficient vector that captures variations in the feature space. The variables $\xi_i$ and $\xi_i^*$ are non-negative slack variables that measure the deviation between the predicted value $\theta_i$ and the observed response $y_i$, from below and above, respectively. The objective function minimizes an asymmetrically weighted sum of these deviations to reflect the targeted quantile level. Constraint (\ref{conv_convex cons}) imposes convexity on the estimated regression function by enforcing a global convexity condition on the local linear approximations.


\section{Convex quantile regression with $\ell_2$-norm regularization}\label{sec:CQR}
Although ridge regression is a widely used method for mitigating overfitting, limited research has explored the impact of ridge regularization on convex quantile regression and its solution methodologies. Therefore, based on the model (\ref{conv_CQR}), we formulate the convex quantile regression with $\ell_2$-norm regularization and $\varepsilon$-insensitive zone as follows:
\begin{subequations}
    \label{CQR}
    \begin{align}
         \min_{\boldsymbol{\beta}, \boldsymbol{\theta} , \xi, \xi^*} & \quad \frac{1}{2} \sum_{i} \|\boldsymbol{\beta}_i\|_2^2 + C \sum_{i=1}^{n} (\tau \xi_i + (1 - \tau) \xi_i^*) \\ \label{positive}
    \text{s.t.} & \quad y_i - \theta_i \leq (1-\tau)\varepsilon + \xi_i  \quad \forall i \in [n], \\ \label{negative}
                & \quad \theta_i - y_i \leq \tau \varepsilon + \xi_i^* \quad \forall i \in [n], \\   \label{convex cons}
                 &\quad \theta_i + \boldsymbol{\beta}_i'(\mathbf{x}_j-\mathbf{x}_i) \leq \theta_j \quad \forall i,j \in [n],\\
                & \quad \xi_i \geq 0, \quad \xi_i^* \geq 0 \quad \forall i \in [n],
    \end{align}
\end{subequations}
where $C$ is a prespecified parameter that controls the trade-off between model complexity and prediction accuracy. This formulation incorporates both $\ell_2$-norm regularization on the subgradients and an $\varepsilon$-insensitive zone ((\ref{positive})--(\ref{negative})), a mechanism from support vector regression that ignores small residuals within a threshold $\varepsilon$ and penalizes only larger deviations via slack variables $\xi_i$ and $\xi_i^*$ \citep{liao2024convex, awad2015support}.

We briefly outline the derivation of the proposed convex quantile regression formulation~(\ref{CQR}). The theoretical and Bayesian motivations for the regularization terms, as well as their connection to Lipschitz convex regression~\citep{mazumder2019computational}, are detailed in Section S.2 of the Supplementary material.

(a) The inclusion of the $\ell_2$-norm regularization serves two main purposes:
\begin{itemize}
    \item To mitigate overfitting by shrinking the local subgradients;
    \item To induce strong convexity in the objective function, which is essential for enabling the Benders decomposition for SCQR in Section 4.
\end{itemize}

(b) The $\varepsilon$-insensitive loss, widely used in support vector regression~\citep{liao2024convex} and quantile regression~\citep{Anand2020}, is introduced here for the first time in convex quantile regression. We empirically assess whether the $\varepsilon$-insensitive zone, originally developed to enhance robustness in support vector regression, can similarly reduce overfitting in convex quantile regression.

Together, the $\ell_2$-norm regularization and the $\varepsilon$-insensitive zone enhance model stability and generalization. When combined with an $\ell_0$-based sparsity constraint, they contribute to more effective variable selection. In Section 4, we introduce the $\ell_0$-constraint and explain its interaction with these regularization components within our convex quantile regression framework.

\subsection{Primal cutting-plane algorithm}
Similar to ordinary convex regression, the $n(n-1)$ convexification constraints~(\ref{convex cons}) substantially increase the computational complexity of the model. To mitigate this, we adapt the cutting-plane algorithm~\citep{balazs2015near, bertsimas2021sparse}, which begins with a small subset of constraints and iteratively adds violated ones in a delayed fashion. At each iteration, we solve a reduced master problem---identical to the full model~(\ref{CQR}) in objective and variables, but with only a subset of constraints. Violated constraints are identified by solving a separation problem for the relaxed solution $(\hat{\boldsymbol{\theta}}, \hat{\boldsymbol{\beta}})$. For each $i \in \{1, \ldots, n\}$, we find
\[
j(i) = \arg \max_{1 \leq k \leq n} \left\{ \hat{\theta}_i + \hat{\boldsymbol{\beta}}_i'(\mathbf{x}_k - \mathbf{x}_i) - \hat{\theta}_k \right\},
\]
and add the corresponding constraint $\theta_i + \boldsymbol{\beta}_i'(\mathbf{x}_{j(i)} - \mathbf{x}_i) \leq \theta_{j(i)}$ to the reduced master problem. The complete procedure is shown in Algorithm~\ref{cutting-plane}.

\begin{algorithm}
\small
\caption{Cutting-Plane Algorithm for Problem (\ref{CQR})}
\label{cutting-plane}
\KwIn{Data \((y_i, \mathbf{x}_i), i = 1, \ldots, n\), tolerance \(Tol > 0\)}
\KwOut{An optimal solution \((\boldsymbol{\hat{\theta}}, \boldsymbol{\hat{\beta}},    \hat{\xi}_1^*, \ldots, \hat{\xi}_n^*, \hat{\xi}_1, \ldots, \hat{\xi}_n)\) to (\ref{CQR})}
\BlankLine
1. Solve the initial reduced master problem \\
2. Set \textit{Continue} \(= True\). \\
3. \While{\textit{Continue} \(== True\)}{
    4. \For{\(1 \leq i \leq n\)}{
        Solve the separation problem and add the corresponding violated constraints to the reduced master problem.\\
    }
    7. \If{there is no violated constraint within the tolerance \(Tol\)}{
        Set \textit{Continue} \(\leftarrow False\). \\
    }
    8. \Else{
        Resolve the updated reduced master problem with  additional constraint(s) added from current iteration \\
    }
 
}
\end{algorithm}

\begin{rmk}
    In practice, the separation step in Line 4 of Algorithm~\ref{cutting-plane}, which involves identifying violated constraints for each sample \( i \in \{1, \ldots, n\} \), can be significantly accelerated using \textbf{parallel computing} or \textbf{matrix-based vectorized operations}.
\end{rmk}
In the next section, we demonstrate that the convex quantile regression problem serves as the Benders subproblem in the GBD algorithm for the SCQR problem. 

\section{Sparse convex quantile regression}\label{sec:Method}
In this section, we introduce the SCQR problem and present the design of the GBD algorithm to address it. The model is reformulated as follows:

\begin{subequations}
\label{SCQR}
\begin{align}
    \min \quad &\frac{1}{2} \sum_{i} \|\boldsymbol{\beta}_i\|_2^2 + C\sum_{i}(\tau \xi_{i} + (1 - \tau)\xi_{i}^*)\\
    \text{s.t.} \quad & y_i - \theta_i \leq  (1-\tau)\varepsilon + \xi_i \quad \forall i \in [n],\\
    &\theta_i - y_i \leq \tau \varepsilon+ \xi_i^* \quad \forall i \in [n],\\ \label{covexcon}
    &\theta_i + \boldsymbol{\beta}_i'(\mathbf{x}_j-\mathbf{x}_i) \leq \theta_j \quad \forall i,j \in [n],\\ \label{Big-M}
    &|(\boldsymbol{\beta}_i)_j| \leq M z_j\quad \forall i \in [n], j \in [d], \\ \label{card}
                &\ \sum_{j=1}^{d} z_j \leq k, \\ \label{z-0-1}
    &\mathbf{z} \in \{0, 1\}^d, \\
                & \xi_i \geq 0, \quad \xi_i^* \geq 0, \quad \forall i \in [n].
\end{align}
\end{subequations}
The SCQR model proposed by \citet{dai2023variable} extends classical CQR (\ref{conv_CQR}) by adding cardinality constraints~(\ref{Big-M})–(\ref{z-0-1}) to perform variable selection, where sparsity is imposed solely through hard constraints. In our work, we build on SCQR by further introducing $\ell_2$-norm regularization and an $\varepsilon$-insensitive zone. Unlike the cardinality constraint, $\ell_2$ regularization does not induce sparsity; instead, it only improves generalization and estimation stability. The combination of these elements not only enhances variable selection accuracy and predictive performance (See the theoretic motivations in Section 3), but also yields the structural properties required to design a tractable decomposition algorithm (Theorem~\ref{min_to_max}), thereby making a significant step toward overcoming SCQR’s computational challenges posed in \cite{dai2023variable}.

 As in conventional (convex) quantile regression, the quantile property in terms of the optimal solution  $\hat{\xi}_i$ and $\hat{\xi}^*_i$ to (\ref{SCQR}) remains essential in this context as well, with its definition provided in \cite{wang2014nonparametric} and \cite{dai2023generalized}. Consequently, the model (\ref{SCQR}) is expected to satisfy the extended quantile property:
\begin{thm}
Let \( n^{-} \) and \( n^{+} \) denote the numbers of observations with strictly negative residuals (i.e., \( \hat{\xi^*_i} > 0 \)) and strictly positive residuals (i.e., \( \hat{\xi}_i > 0 \)), respectively. Then, the following quantile property holds:
\begin{equation}
\frac{n_-}{n} \leq \tau \quad \text{and} \quad \frac{n_+}{n} \leq 1 - \tau.
\end{equation}
\end{thm}
The proofs of this theorem and others may be found in Section S.1 of the Supplementary material.

\begin{figure}[htbp]
	\centering
	\begin{minipage}{0.5\linewidth}
		\vspace{3pt}
		\centerline{\includegraphics[width=\textwidth]{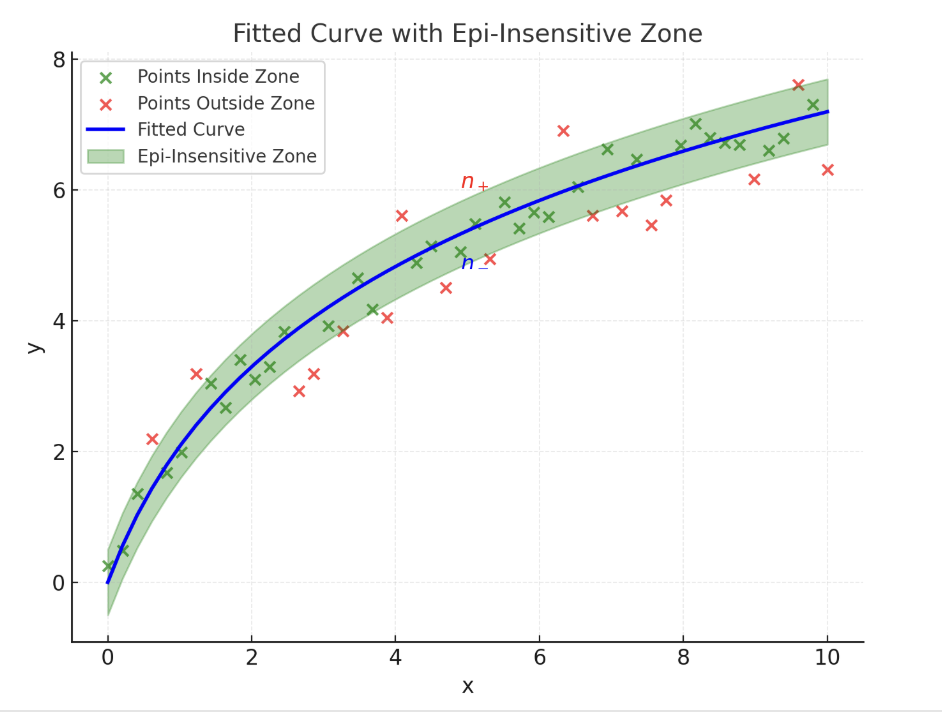}}
	\end{minipage}
	\caption{$\varepsilon$-insensitive zone and illustration of quantile property}
	\label{fig:example}
\end{figure}

\begin{rmk}
The quantile property in this theorem differs from the classical formulation due to the presence of the $\varepsilon$-insensitive zone. In standard convex quantile regression, residuals directly determine the quantile property. In contrast, our model defines residuals through the slack variables $\xi_i$ and $\xi_i^*$, which are strictly positive only when the prediction error lies outside the $\varepsilon$-insensitive zone $[-\tau \varepsilon, (1-\tau)\varepsilon]$. Thus, the quantile property here describes the proportion of observations with nonzero slack variables—i.e., those whose residuals exceed the tolerance range. This reformulation reflects how the $\varepsilon$-insensitive region modifies the classical residual distribution. Figure~\ref{fig:example} illustrates this mechanism.
\end{rmk}

\begin{rmk}
    For clarity of presentation, we omit the consideration of the $\varepsilon$-insensitive zone in the following discussion, as it has no impact on the derivation of our theoretical results or algorithms.
\end{rmk}

The model’s complexity stems from three main sources. 
First, constraints such as the cardinality constraint~(\ref{card}) and Big-M constraints~(\ref{Big-M}) substantially increase computational burden. 
Second, the convexification constraints~(\ref{covexcon}) add numerous restrictions to enforce convexity, further complicating the formulation. 
Finally, selecting an appropriate Big-M constant~(\(M\)) is challenging: overly large values over-relax the model, while overly small values risk excluding feasible solutions~\citep{bertsimas2016best}. 
These factors together make the model design and implementation intricate.


\subsection{Generalized Benders decomposition}
By fixing $\mathbf{z}$, the problem reduces to a standard CQR model (\ref{CQR}), which naturally separates the combinatorial task of variable selection from the functional fitting of CQR. This reformulation avoids tackling both sources of difficulty simultaneously in model~(\ref{SCQR}) and instead casts them into a master problem and a subproblem that can be solved more efficiently. The Benders decomposition framework is particularly suitable in this setting, as it iteratively coordinates the two problems through generated cuts, thereby enhancing computational efficiency and scalability compared with solving the original mixed-integer formulation directly. For a given $k$, we introduce $S_k^d$ as the set of d-dimensional binary vectors with at most k nonzero components; that is: $S_k^d = \left\{ \mathbf{z} \in \{0,1\}^d : \sum_{i=1}^{d} z_i \leq k \right\}$.

Then the model (\ref{SCQR}) can be reformulated as:
\begin{equation} \label{formula1}
    \min_{\mathbf{z} \in S_k^d} \quad g(\mathbf{z}), 
\end{equation}
where 
\begin{subequations}
\label{subproblem}
\begin{align}
    g(\mathbf{z}) = \min_{\boldsymbol{\beta}_i, \mathbf{\theta}, \xi, \xi^*, \mathbf{Z} = diag(\mathbf{z})} & \quad \frac{1}{2} \sum_{i=1}^{n} \|\boldsymbol{\beta}_i\|_2^2 + C \sum_{i=1}^{n} (\tau \xi_i + (1 - \tau)\xi_i^*) \\
    \text{s.t.}& \quad y_i - \theta_i \leq \xi_i \quad \forall i \in [n], \\
                & \quad \theta_i - y_i \leq \xi_i^* \quad \forall i \in [n], \\ \label{dual_convex}
                & \quad \theta_i + \boldsymbol{\beta}_i'\mathbf{Z}(\mathbf{x}_j- \mathbf{x}_i) \leq \theta_j\quad \forall i,j \in [n],\\
                & \quad \xi_i \geq 0, \quad \xi_i^* \geq 0 \quad \forall i \in [n].
\end{align}
\end{subequations}
The next theorem shows that the minimization subproblem $g(\mathbf{z})$, which involves the piecewise-linear, non-smooth, and asymmetric structures, can be reformulated as a maximization problem. This reformulation enables the derivation of the critical Benders cuts.

\begin{thm}\label{min_to_max}
    The problem (\ref{formula1}) is equivalent to solving the following formulation with binary variables and convex objective.
    \begin{equation}\label{max-formula}
        \min_{\mathbf{z} \in S_k^d} \quad g(\mathbf{z}), 
    \end{equation}
    where
    \begin{subequations}
    \label{inner}
    \begin{align}
    g(\mathbf{z}) = \max \quad &-\frac{1}{2} \sum_{i} \|\sum_{j}\mu_{ij}\textbf{Z}(\mathbf{x}_j - \mathbf{x}_i)\|_2^2 + \sum_{i}\lambda_i y_i - \sum_{i}\lambda_i^* y_i\\
    \text{s.t.} \quad & -\lambda_i + \lambda_i^* + \sum_{j}\mu_{ij} - \sum_{j}\mu_{ji} = 0 \quad \forall i \in [n],\\
    &0\leq \lambda_i \leq \tau C \quad \forall i \in [n],\\
    &0\leq \lambda_i^* \leq (1-\tau) C \quad \forall i \in [n],\\
    &\mu_{ij} \geq 0  \quad \forall i,j \in [n].
    \end{align}
 \end{subequations}
    
\end{thm}

According to the theorem, \(g(\mathbf{z})\) is a convex function with its subgradient $\partial g(\mathbf{z})$ at point \(\mathbf{z}\) given by \(-\frac{1}{2}\sum_{i=1}^n \left( \sum_{j=1}^n \hat{\mu}_{ij} (\mathbf{x}_j - \mathbf{x}_i) \right)^2\), where $\hat{\mu}$ is the optimal dual solution to (\ref{inner}). Consequently, we can reformulate (\ref{formula1}) into the Benders formulation based on this fact.

\begin{thm}
    The formulation (\ref{formula1}) can be transformed into the Benders formulation:
    \begin{subequations}
    \label{benders formulation}
    \begin{align}
\min_{\mathbf{z} \in \{0,1\}^d, \gamma} & \quad \gamma \\ \label{benders cut}
\text{s.t.} & \quad g(\mathbf{z}^*) + \partial g(\mathbf{z}^*)'(\mathbf{z} - \mathbf{z}^{*}) \leq \gamma \quad \forall \mathbf{z}^* \in S^d_k, \\
& \quad \sum_{i=1}^d z_i \leq k,
\end{align}
\end{subequations}
    
\end{thm}

In Benders decomposition, constraints (\ref{benders cut}) are called the Benders cuts and  problem (\ref{benders formulation}) is solved using the delayed constraint generation algorithm. The full model of problem (\ref{benders formulation}) is referred to as the master problem, while the model containing only a subset of the constraints in (\ref{benders cut}) is known as the reduced master problem.We begin by solving the initial reduced master problem. Next, we identify the violated Benders cuts by solving the Benders subproblem (\ref{subproblem}) and iteratively incorporate them in a delayed manner. At each iteration, the updated reduced master problem is solved with the newly added violated Benders cuts, gradually refining the solution. The reduced master problem at iteration $t$ can be formulated as follows:\
\begin{subequations}
    \label{reduced master problem}
    \begin{align}
\min_{\mathbf{z} \in \{0,1\}^d, \gamma} & \quad \gamma \\
\text{s.t.} & \quad g(\mathbf{z}^*) + \partial g(\mathbf{z}^*)'(\mathbf{z} - \mathbf{z}^{*}) \leq \gamma \quad \forall \mathbf{z}^* \in S^t, \\
& \quad \sum_{i=1}^d z_i \leq k,
\end{align}
\end{subequations}
where $S^t$ denotes the collection of all feasible solutions identified up to iteration $t$. It is worth noting that the Benders subproblem (\ref{subproblem}) corresponds precisely to the CQR model (\ref{CQR}). Therefore, the cutting-plane algorithm \ref{cutting-plane} serves as an effective tool for solving the Benders subproblem. Using this approach, we can obtain the optimal values of the dual variables \(\mu\) associated with the constraints (\ref{dual_convex}) and subsequently compute the required gradients.

\subsection{Warm start approach}
As mentioned in the literature \citep{bertsimas2021sparse, bertsimas2020sparse}, in the context of sparse regression for standard mean value regression, the linear relaxation of (\ref{max-formula}) offers a relatively tight approximation to problem (\ref{SCQR}) and therefore may provide good-quality warm starts. Building on this result, we extend the conclusion to our context and derive the following corollary.

\begin{corollary}
    The linear relaxation of (\ref{max-formula}) can be characterized by the following optimization problem:
    \begin{subequations}
    \label{linear relaxation}
    \begin{align}
\min_{\boldsymbol{\mu} \geq 0, \lambda, \lambda^*, \gamma } & \quad \gamma- \sum_{i}\lambda_i y_i + \sum_{i}\lambda_i^* y_i \\ 
\text{s.t.} & \quad \gamma \geq \frac{1}{2}\sum_{p = 1}^d z_p\left\{ \sum_{i=1}^n \left( \sum_{j=1}^n \mu_{ij} (\mathbf{x}_j - \mathbf{x}_i) \right)^2_p \right\} \forall \mathbf{z} \in S^d_k,\\
 &\quad  -\lambda_i + \lambda_i^* + \sum_{j}\mu_{ij} - \sum_{j}\mu_{ji} = 0 \quad \forall i \in [n],\\
 &\quad 0\leq \lambda_i \leq \tau C \quad \forall i \in [n],\\
&\quad 0\leq \lambda_i^* \leq (1-\tau) C \quad \forall i \in [n],\\
&\quad \mu_{ij} \geq 0 \quad \forall i,j \in [n],
\end{align}
\end{subequations}
\end{corollary}
This problem can also be solved through the cutting-plane algorithm. And the initial support set of $\mathbf{z}$ would be the corresponding indices of the largest $k$ values of the components of vector  $\sum_{i=1}^n \left( \sum_{j=1}^n \hat{\mu}_{ij} (\mathbf{x}_j - \mathbf{x}_i) \right)^2$. Now we can display the whole GBD algorithm here:

\begin{algorithm}
\caption{Generalized Benders Decomposition}
\label{Bender algorithm}
\KwIn{$C > 0$, $T > 0$}
\KwOut{Optimal support $\mathbf{z}^*$, lower bound $LB$, upper bound $UB$}
\BlankLine
1. Start with $\gamma^0 = 0$, initial feasible $\mathbf{z}^0$ via warm start\;
2. Set $t \gets 0$, initialize $LB \gets -\infty$, $UB \gets +\infty$, $\mathbf{z}^* \gets \mathbf{z}^0$\;
3. \While{$UB - LB > 0$ \textbf{and} $t \leq T$}{
    4. Solve subproblem at $\mathbf{z}^t$ to compute $g(\mathbf{z}^t)$ and subgradient $\partial g(\mathbf{z}^t)$ via Theorem \ref{min_to_max}\;
    5. \textbf{Update upper bound:} $UB \gets \min(UB, g(\mathbf{z}^t))$\;
    6. Add cutting-plane constraint: $g(\mathbf{z}^t) + \partial g(\mathbf{z}^t)'(\mathbf{z} - \mathbf{z}^t) \leq \gamma$\;
    7. Resolve the reduced master problem (\ref{reduced master problem}) to obtain $(\mathbf{z}^{t+1}, \gamma^{t+1})$\;
    8. \textbf{Update lower bound:} $LB \gets \max(LB, \gamma^{t+1})$\;
    9. \textbf{Update incumbent:} If $g(\mathbf{z}^{t}) < g(\mathbf{z}^*)$, set $\mathbf{z}^* \gets \mathbf{z}^{t}$\;
    10. $t \gets t + 1$\;
}
\end{algorithm}

\begin{rmk}
   Although our primary focus is on convex quantile regression, the proposed algorithm has broader applicability. In particular, it can be directly extended to address sparse linear quantile regression as a special case, highlighting its versatility in handling a wider range of high-dimensional quantile modeling tasks. Furthermore, the GBD algorithm specifically designed for our problem is an exact decomposition method that converges to the optimal solution in a finite number of iterations.
\end{rmk}

\begin{thm}
    Algorithm \ref{Bender algorithm} can converge to the optimal solution in a finite number of iterations.
\end{thm}

\begin{rmk}
While the finite convergence of our algorithm is theoretically guaranteed, establishing a general convergence rate remains challenging due to the NP-hard nature of the SCQR problem \citep{dai2023variable}. In the worst case, the algorithm may require an exponential number of iterations to reach optimality, as commonly encountered in integer programming \citep{wolsey1999integer, rahmaniani2017benders}. Nonetheless, our computational results show that the algorithm performs efficiently in practice and consistently yields high-quality solutions across a range of instances.
\end{rmk}

\section{Local search-based Benders matheuristic}\label{heuristic}
In our preliminary experiments, we observed that although the GBD algorithm may struggle to tighten the lower bound and reach convergence on larger problems, it consistently finds high-quality solutions in just a few iterations. This makes GBD a promising matheuristic for real-world applications.

To further refine the incumbent solution, we propose a novel \textit{local search-based Benders} (LSB) matheuristic. Local search is a well-established strategy for combinatorial optimization \citep{lourencco2003iterated}, and recent work has explored its integration with Benders decomposition to accelerate cut generation or global convergence \citep{rei2009accelerating, maher2021enhancing, fischetti2003local}. However, these approaches do not fully leverage Benders decomposition to explore solution neighborhoods directly. Our LSB method addresses this gap by using Benders decomposition as the engine to solve localized subproblems, making it both simple and effective in improving the incumbent. To the best of our knowledge, such a combination has not been explicitly studied in the literature.

Given the incumbent integer solution \( \mathbf{z}^* \) obtained from Algorithm \ref{Bender algorithm}, we define the neighborhood around \( \mathbf{z}^* \) within a predefined distance \( r \) as
\[
\mathcal{N}(\mathbf{z}^*, r) = \{\mathbf{z} \in \{0,1\}^d : d(\mathbf{z}, \mathbf{z}^*) \leq r\},
\]
where \( d(\mathbf{z}, \mathbf{z}^*) \) represents the Hamming distance \citep{bookstein2002generalized} between the two binary vectors \( \mathbf{z} \) and \( \mathbf{z}^* \). The Hamming distance, \( d(\mathbf{z}, \mathbf{z}^*) \), is defined as the number of positions at which the corresponding bits of \( \mathbf{z} \) and \( \mathbf{z}^* \) differ. This neighborhood forms the feasible search space for exploring alternative solutions close to \( \mathbf{z}^*\). Then the restricted master problem is as follows:
\begin{subequations}
    \begin{align}
        \min_{\mathbf{z} \in \{0,1\}^d, \gamma} & \quad \gamma \\ \label{res benders cut}
\text{s.t.} & \quad g(z^*) + \partial g(\mathbf{z}^*)'(\mathbf{z} - \mathbf{z}^{*}) \leq \gamma, \quad \forall \mathbf{z}^* \in S^d_k, \\
& \quad \sum_{i=1}^d z_i \leq k,\\
& \quad \mathbf{z} \in \mathcal{N}(\mathbf{z}^*, r).
    \end{align}
\end{subequations}
Similarly, we refer to the problem with partial constraints in (\ref{res benders cut}) as the reduced restricted master problem. Here we present the complete algorithm of the LSB approach in Algorithm \ref{LSB algorithm}.
\begin{algorithm}
\small
\caption{Local Search-Based Benders Matheuristic}
\label{LSB algorithm}
\KwIn{Incumbent solution \(\mathbf{z}^*\), Hamming distance \(r > 0\), maximum iterations \(T > 0\)}
\KwOut{Improved incumbent solution \(\mathbf{z}^*\)}
\BlankLine
1. Define the neighborhood \(N(\mathbf{z}^*, r) = \{\mathbf{z} \in \{0, 1\}^d : d(\mathbf{z}, \mathbf{z}^*) \leq r\}\), where \(d(\mathbf{z}, \mathbf{z}^*)\) is the Hamming distance. \\
2. Set \(t \gets 0\). \\
3. \While{Termination criteria are not met and \(t \leq T\)}{
    4. Apply Benders decomposition (Algorithm \ref{Bender algorithm}) within the neighborhood \(N(\mathbf{z}^*, r)\) to explore feasible solutions. \\
    5. Update the solution \(\mathbf{z}^*\) if a better solution is found. \\
    6. Redefine the neighborhood \(N(\mathbf{z}^*, r)\) based on the updated \(\mathbf{z}^*\). \\
    7. \(t \gets t + 1\). \\
}
\end{algorithm}

The proposed LSB matheuristic integrates the simplicity of local search with the decomposition power of Benders methods to effectively explore feasible solutions. A key innovation lies in introducing a localized constraint, which mitigates one of the major challenges in Benders decomposition—oscillations caused by excessive exploration of distant, suboptimal regions \citep{rahmaniani2017benders}. By restricting the search to a targeted neighborhood, the method stabilizes the solution process and improves computational efficiency. This synergy offers a promising heuristic framework for solving large-scale combinatorial optimization problems, while also creating opportunities to incorporate enhancement techniques traditionally developed for local search.

\section{Simulation study}\label{sec:exp}






In all the experiments that follow, we use Gurobi 10.0.3 as the optimization solver, running on a MacBook Pro 14-inch (2021) equipped with an Apple M1 Pro chip and 16 GB memory, under macOS Sonoma 14.0. The experiments are implemented using the Python programming language with the Gurobipy interface for model formulation and solution. Table \ref{tab:notation} summarizes the parameters and notations introduced in this section.

\begin{table}[h!]
\centering
\small
\caption{Notation and Descriptions}
\label{tab:notation}
\begin{tabular}{clcl}
\hline
\textbf{Sym.} & \textbf{Description} & \textbf{Sym.} & \textbf{Description} \\
\hline
$\mathbf{X}$    & Feature matrix        & $n$             & Number of data points         \\
$d$             & Total number of features & $\varepsilon$  & $\varepsilon$-insensitive zone parameter \\
$\tau$          & Quantile level (e.g., 0.5 for median) & $\rho$ & Feature correlation \\
$C$             & Penalty coefficient   & $k$             & Number of selected features \\
\hline
\end{tabular}
\end{table}

\subsection{Test for formulation (\ref{CQR}) in reducing overfitting}
In this subsection, we consider two data generating processes (DGP) (see, e.g., \cite{liao2024convex} ):
\begin{align*}
&\text{(1) DGP I:} \quad  y = 3 + x_1^{0.2} + x_2^{0.3} + \epsilon, \\
&\text{(2) DGP II:} \quad  y = 3 + x_1^{0.05} + x_2^{0.15} + x_3^{0.3} + \epsilon,
\end{align*}
where \( x_1, x_2, x_3\) are independently and randomly sampled from the uniform distribution \( U[1, 10] \) and the error term \( \epsilon \) is drawn from \( \mathcal{N}(0, \sigma^2) \). For each DGP, we consider  different scenarios with $n \in \{100, 500\}$ and $\sigma = 1$. For each scenario, we replicate 10 times to calculate the in-sample and out-of-sample Mean Absolute Error (MAE). In the context of quantile regression, the out-of-sample MAE on a test set \( \mathcal{D}_{\text{test}} = \{ (x_i^{\text{test}}, y_i^{\text{test}}) \}_{i=1}^{n_{\text{test}}} \) is defined as:

\[
\text{MAE} = \frac{1}{n_{\text{test}}} \sum_{i=1}^{n_{\text{test}}} \rho_\tau \left( y_i^{\text{test}} - \hat{y}_i^{\text{test}} \right),
\]
where \( \hat{y}_i^{\text{test}} \) denotes the predicted \(\tau\)-quantile of the conditional distribution of \( y \) given \( x_i^{\text{test}} \), while the in-sample MAE is similarly computed using the quantile loss function on the training set. We select the model parameters \( C \) and \( \varepsilon \) using five-fold cross-validation, where \( C \) and \( \varepsilon \) are chosen from the sets \( \{0.1, 0.5, 1, 2, 5\} \) and \( \{0, 0.02, 0.2, 1, 2\} \), respectively.

To assess the roles of the $\ell_2$-norm regularization and the $\varepsilon$-insensitive zone in mitigating overfitting, we first evaluate a variant with $\varepsilon=0$ and only the $\ell_2$ term (\textbf{CQR-$\ell_2$}). We then introduce $\varepsilon$ to examine their combined effect (\textbf{CQR-$\ell_2$-$\varepsilon$}). For comparison, we also consider the Lipschitz convex quantile regression (\textbf{LCQR}) from \citet{mazumder2019computational}, where Lipschitz constraints are applied directly to convex quantile regression (see Section S.2.1 Formulation~(S2) in the Supplementary material), and the baseline convex quantile regression (\textbf{CQR}) without regularization. These comparisons highlight the effectiveness of our techniques in reducing overfitting.

\begin{table}[ht]
\centering
\small
\caption{In-sample (In) and Out-of-sample (Out) MAE comparison with $\sigma = 1$ and $\tau = 0.25$.}
\label{tab:combined_mae}
\begin{tabular}{
cc
cc
cc
cc
cc}
\toprule
\textbf{DGP} & \textbf{$n$} 
& \multicolumn{2}{c}{\textbf{CQR}} 
& \multicolumn{2}{c}{\textbf{CQR-$\ell_2$}} 
& \multicolumn{2}{c}{\textbf{CQR-$\ell_2$-$\varepsilon$}} 
& \multicolumn{2}{c}{\textbf{LCQR}} \\
 & & In & Out & In & Out & In & Out & In & Out \\
\midrule
I   & 100 & 0.270 & 0.455 & 0.298 & 0.323 & 0.305 & \textbf{0.320} & 0.294 & 0.324 \\
    & 500 & 0.300 & 0.341 & 0.311 & 0.314 & 0.310 & \textbf{0.312} & 0.312 & 0.314 \\
II  & 100 & 0.206 & 0.658 & 0.290 & 0.323 & 0.299 & \textbf{0.322} & 0.294 & \textbf{0.322} \\
    & 500 & 0.274 & 0.670 & 0.309 & 0.317 & 0.309 & \textbf{0.316} & 0.314 & \textbf{0.316} \\
\bottomrule
\end{tabular}
\end{table}

Table~ \ref{tab:combined_mae} reports the in-sample and out-of-sample MAE at $\tau = 0.25$, averaged over ten trials. The results show that standard CQR suffers from overfitting, with higher out-of-sample MAE than regularized variants. Adding $\ell_2$ regularization (CQR-$\ell_2$) markedly improves performance, and incorporating the $\varepsilon$-insensitive zone (CQR-$\ell_2$-$\varepsilon$) further stabilizes results. Compared to LCQR, our method attains similar or better out-of-sample accuracy without additional Lipschitz constraints, and retains the structural properties necessary for our decomposition-based algorithm.



\subsection{Monte Carlo study related to the GBD algorithm}
In this subsection, we present numerical experiments to evaluate the performance of the core algorithm proposed in this paper, namely the GBD method. The experiments are designed to examine two main aspects: (1) the computational efficiency of solving the Benders subproblem (i.e, problem (\ref{CQR})), and (2) the overall effectiveness and accuracy of the full algorithm in performing variable selection.

\subsubsection{Data description}\label{sec:data}

We generate the synthetic data for our next experiments using the following procedure (see, .e.g, \cite{bertsimas2021sparse}). The feature matrix $\mathbf{X}$ is generated from a standard Gaussian distribution. The response variable $y_i$ is modeled using the convex function $\Phi(\mathbf{x}) = \|\mathbf{x}\|_2^2$, with an additive Gaussian noise $\epsilon_i$, defined as: $y_i = \Phi(\mathbf{x}_i) + \epsilon_i, \epsilon_i \sim N(0, \sigma^2), i = 1, \ldots, n.$
Here, the errors $\epsilon_i$ are assumed to be independent and identically distributed (i.i.d.). The variance of the noise $\sigma^2$ is determined by the signal-to-noise ratio (SNR), defined as: $\text{SNR} = \frac{\mathrm{Var}(\mu)}{\mathrm{Var}(\epsilon)}$, where $\mu_i = \Phi(\mathbf{x}_i)$. A higher SNR indicates smaller noise levels relative to the signal, leading to less distortion in the observed data.\

We will report the number of cuts added at each iteration when implementing Algorithm \ref{cutting-plane} and the metric called primal infeasibility \citep{mazumder2019computational}: $\text{Primal infeasibility} = \frac{1}{n} \| \mathbf{V} \|_F$,
where the matrix $\mathbf{V}$ is defined with entries $
    V_{ij} = \max\{0, \hat{\theta}_i + \hat{\boldsymbol{\beta}}_i^\top (\mathbf{x}_j - \mathbf{x}_i) - \hat{\theta}_j\}, \quad \forall i, j \in \{1, \dots, n\}$.
Here, $V_{ij}$ quantifies the degree of violation of the corresponding constraint, with $V_{ij} = 0$ indicating no violation. The Frobenius norm $\|\cdot\|_F$ is given by: $\| \mathbf{V} \|_F^2 = \sum_{i=1}^n \sum_{j=1}^n V_{ij}^2$.

\subsubsection{Test for  Algorithm \ref{cutting-plane}}\label{sec:algorithm1}
In this section, we report the running time and primal infeasibility of Algorithm \ref{cutting-plane} when solving the convex quantile regression with $\ell_2$-norm regularization (\ref{CQR}).  To thoroughly evaluate the performance of Algorithm \ref{cutting-plane}, we conduct experiments across varying parameters, including the quantile value ($\tau$), the size of the training dataset ($n$) and the number of features ($d$) used in Algorithm \ref{cutting-plane}. 

The training dataset sizes $n$ are chosen from \{2000, 10000, 20000\}, while the number of features $d$ is chosen from  the set \{50, 70, 90\}. The tolerance parameter $Tol$ is set to be 0.01, the SNR is set to be 3 and  the quantile level $\tau$ is tested at \{0.25, 0.5, 0.75\}. We set the regularization parameter $C$ to 10. For each combination of these parameters, we record the time required for Algorithm \ref{cutting-plane} to converge and calculate the corresponding primal infeasibility. The reported results represent the average of five independent runs, each using randomly generated data.

\begin{table}[h!]
\centering
\caption{Run time and primal infeasibility of convex quantile regression}
\label{tab:primal_inf}
\scriptsize
\begin{tabular}{ccccccc}
\toprule
$\tau$ & $n$ & $d$ & iter & Infeasibility & Run time (s) \\
\midrule
0.25 & 2000 & 50 & 36 & 0.0008 & 157 \\
     & 10000 & 70 &49 & 0.0010 & 1689 \\
     & 20000 & 90 &  59& 0.0230 & 5551 \\
0.50 & 2000 & 50 & 35 &0.0008& 129 \\
     & 10000 & 70 & 52 & 0.0151 &1605  \\
     & 20000 & 90 & 62 & 0.0362 &  6361 \\
0.75 & 2000 & 50 &38  & 0.0018 &  169 \\
     & 10000 & 70 & 54 & 0.0301 & 2365 \\
     & 20000 & 90 & 68 & 0.0768 & 8366 \\
\bottomrule
\end{tabular}
\end{table}


The results are presented in Table~\ref{tab:primal_inf} and several observations can be obtained:
\begin{itemize}

\item \textbf{Impact of Problem Size:} As $n$ and $d$ increase, the problem becomes more computationally demanding due to more iterations and constraints.

\item \textbf{Algorithm Efficiency:} The proposed cutting-plane algorithm consistently solves all tested instances within minutes, demonstrating strong scalability.


\item \textbf{Effect of Quantile Level:} Larger quantile levels $\tau$ lead to increased run times, suggesting added optimization complexity.
\end{itemize}

Furthermore, we present additional visualizations to provide insights into the iterative behavior of Algorithm \ref{cutting-plane} in Section S.3.1, Figure S1 of the Supplementary material.

\subsubsection{Test for sparse convex quantile regression}\label{sec:sparse}
In this section, we present the computational results for solving the SCQR problem using our proposed GBD algorithm, enhanced by the LSB algorithm. Specifically, we first run Algorithm~\ref{Bender algorithm} for 80 iterations to obtain the incumbent solution \( \mathbf{z}^* \), which is then passed to Algorithm~\ref{LSB algorithm}. The latter is executed with the parameter \( r \) alternating between 1 and 2 every 30 iterations to balance exploration and exploitation, with a maximum iteration limit of \( T = 300 \).  

To generate the simulation dataset, we sample the matrix \( \mathbf{X} \) from a Gaussian distribution. A support set of size \( k \) is randomly chosen from set \( \{1, \dots, d\} \). For each observation \( i \), \( \mathbf{x}_i \) is drawn from a Gaussian distribution with zero mean and a correlation matrix \( \Sigma \), where the entries are defined as $\Sigma_{ij} = \rho^{|i-j|}, \quad \text{for } 1 \leq i, j \leq d$, with \( \rho \in [0, 1] \) controlling feature correlation. Higher \( \rho \) values indicate stronger correlations among features. To enhance numerical stability and improve prediction accuracy, we mean-center and normalize the features and response vectors to ensure a unit \( \ell_2 \) norm. For model selection, we use cross-validation to choose \( C \) from \( \{0.1, 1, 10, 100\} \). To examine the impact of the \( \varepsilon \)-insensitive zone, we conduct experiments with \( \varepsilon \) values of \( \{0, 0.04\} \). In practice, the optimal \( \varepsilon \) should also be determined through cross-validation or other model selection techniques.  

We evaluate the final solution accuracy as a function of SNR, \( \tau \), \( \rho \), \( \varepsilon \), \( d \), and \( k \). Accuracy is defined as:  
\begin{equation}
    \text{Accuracy} = \frac{|S^* \cap \hat{S}|}{k}
    \label{eq:accuracy}
\end{equation}
where \( S^* \) denotes the true support set, and \( \hat{S} \) represents the estimated optimal set obtained by our algorithm.  

Table \ref{tab:all_results} presents the results for \( n = 800 \) with \( \tau = 0.25 \) (Results for \( 0.5 \), and \( 0.75 \) are presented in Section S.3.2, Table S1 of the Supplementary material.). For each \( \tau \), results are provided for \( SNR = 3 \) and \( SNR = 1 \). We generate synthetic data for each parameter combination, creating five datasets per setting. The reported results represent the average over these five experiments. The \textbf{run time} refers to the total computational time (in seconds) measured until the last heuristic solution update, representing the time required to obtain the final solution.
\setlength{\tabcolsep}{3pt}
\renewcommand{\arraystretch}{0.6}

\begin{table}[h!]
\centering
\scriptsize
\caption{Accuracy, iterations and run time for SCQR (n=800)}
\label{tab:all_results}

\begin{subtable}[t]{0.49\textwidth}
\centering
\caption{$\tau=0.25$, SNR=3}
\label{subtab:025_3}
\resizebox{\textwidth}{!}{
\begin{tabular}{ccccccc}
\toprule
\textbf{$\rho$} & $\varepsilon$ & $d$ & $k$ & accuracy (\%) & iteration & run time (s) \\
\midrule
0.2 & 0 & 100 & 10 &90 &  102 &  1358\\
     &   &     & 20 & 88 & 85  & 1789\\
     &   & 40 & 5 & 87 & 149&1920 \\
     &   &     & 10 & 83 & 251 &2753 \\
 & 0.04 & 100 & 10 & 100 & 83  &  994 \\
     &   &     & 20 & 98 & 110 &  2452\\
     &   & 40 & 5 & 93 & 182 & 2210\\
     &   &     & 10 & 100 & 99&1345 \\
\midrule
0.5 & 0 & 100 & 10 &  92 &  164 & 1932 \\
     &   &     & 20 & 91 &  98 & 2090 \\
     &   & 40 & 5 & 96 & 151 &1814 \\
     &   &     & 10 & 88 & 151&1946 \\
 & 0.04 & 100 & 10 & 96 &119  &1403 \\
     &   &     & 20 & 92 & 134 & 2570\\
     &   & 40 & 5 & 96 & 122&1356 \\
     &   &     & 10 & 96 & 181 & 2404\\
\bottomrule
\end{tabular}%
}
\end{subtable}
\hfill
\begin{subtable}[t]{0.49\textwidth}
\centering
\caption{$\tau=0.25$, SNR=1}
\label{subtab:025_1}
\resizebox{\textwidth}{!}{%
\begin{tabular}{ccccccc}
\toprule
\textbf{$\rho$} & $\varepsilon$ & $d$ & $k$ & accuracy (\%) & iteration & run time (s) \\
\midrule
0.2 & 0 & 100 & 10 &90 & 110 &  1350 \\
     &   &     & 20 &  88& 156 &  3043\\
     &   & 40 & 5 &92 & 144 &1802 \\
     &   &     & 10 & 88 & 105&1204 \\
 & 0.04 & 100 & 10 & 92 & 102 & 1201  \\
     &   &     & 20 & 95 &  178 &  3422\\
     &   & 40 & 5 & 96 & 94&1137 \\
     &   &     & 10 & 100 & 104&1405 \\
\midrule
0.5 & 0 & 100 & 10 & 88 & 97 & 1154 \\
     &   &     & 20 & 86 & 101 & 2056 \\
     &   & 40 & 5 & 88& 197& 2305 \\
     &   &     & 10 & 88 & 166&2411 \\
 & 0.04 & 100 & 10 & 92 & 105 & 1258\\
     &   &     & 20 & 86 & 136 & 2804\\
     &   & 40 & 5 & 92 & 121 &1352 \\
     &   &     & 10 & 92 & 72 & 830 \\
\bottomrule
\end{tabular}%
}
\end{subtable}
\end{table}

From these tables, we can make the following observations:
\begin{itemize}
    \item \textbf{High accuracy:} Our algorithm achieves near 90\% feature selection accuracy, even under low signal-to-noise ratios and high feature correlation.
    \item \textbf{Efficiency and scalability:} The proposed decomposition framework, together with LSB method, rapidly identifies high-quality solutions within a few iterations.
    \item \textbf{Quantile and $\varepsilon$-zone effects:} The  $\varepsilon$-insensitive zone enhances both computational efficiency and estimation accuracy, particularly at lower quantile levels.
\end{itemize}

\subsubsection{Comparision with CNLS-A algorithm}

\begin{figure}[htbp]
\centering
\begin{subfigure}[t]{0.4\textwidth}
  \centering
  \includegraphics[width=\linewidth]{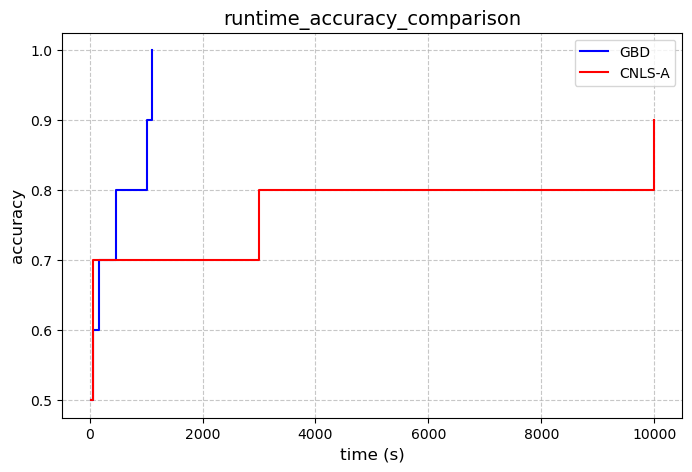}
  \caption{$d = 40$, $k = 10$, SNR=3,  $\rho = 0.2$}
  \label{fig:tau_0.25_40_10_3_0.2}
\end{subfigure}
\hfill
\begin{subfigure}[t]{0.4\textwidth}
  \centering
  \includegraphics[width=\linewidth]{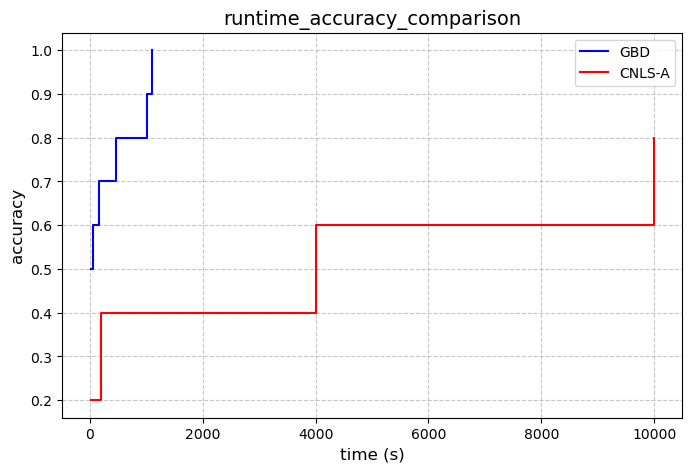}
  \caption{$d = 40$, $k = 10$, SNR=1,  $\rho = 0.5$}
  \label{fig:tau_0.25_40_10_1_0.5}
\end{subfigure}

\vspace{0.2cm} 
\begin{subfigure}[t]{0.4\textwidth}
  \centering
  \includegraphics[width=\linewidth]{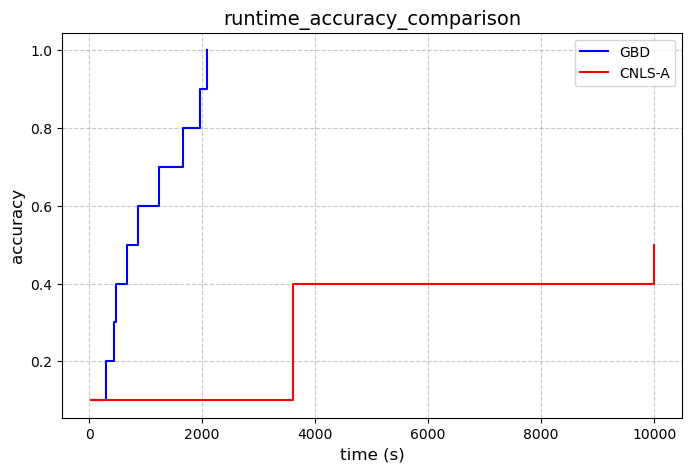}
  \caption{$d = 100$, $k = 10$, SNR=3,  $\rho = 0.2$}
  \label{fig:tau_0.25_100_10_3_0.2}
\end{subfigure}
\hfill
\begin{subfigure}[t]{0.4\textwidth}
  \centering
  \includegraphics[width=\linewidth]{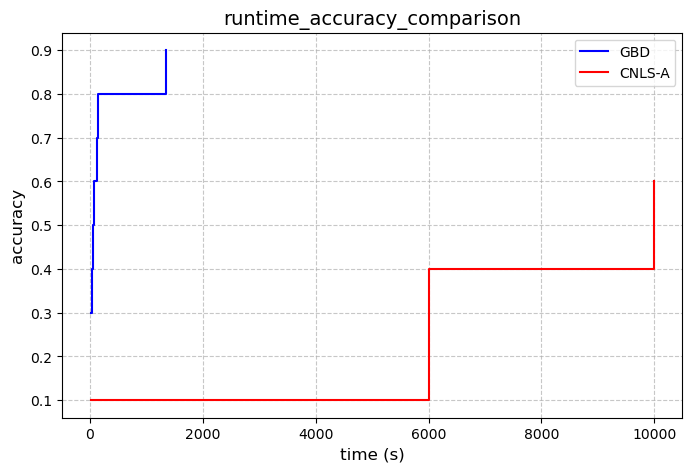}
  \caption{$d = 100$, $k = 10$, SNR=1,  $\rho = 0.5$}
  \label{fig:tau_0.25_100_10_1_0.5}
\end{subfigure}

\caption{Accuracy-time curves for GBD and CNLS-A at $\tau = 0.25$.}
\label{fig:plot_compare}
\end{figure}

A primal cutting-plane algorithm for convex mean regression was proposed by \cite{bertsimas2021sparse}, and was subsequently adapted for convex quantile regression as the CNLS-A algorithm (Algorithm 1 in \cite{dai2023variable}). The CNLS-A algorithm iteratively generates convexity constraints (\ref{covexcon}) in the primal formulation (\ref{SCQR}). This approach is relatively time-consuming, as it requires solving a relaxed version of model (\ref{SCQR}) containing only a subset of the convexity constraints at each iteration \citep{bertsimas2021sparse}. We compare our GBD algorithm with CNLS-A by plotting the evolution of variable selection accuracy over time. For CNLS-A,  we select $M$ over $\{0.1, 1, 5, 10\}$  and set a time limit of \textbf{10,000} seconds and record the incumbent solution accuracy at each time point. Figure~\ref{fig:plot_compare} shows the results for \(\tau = 0.25\), while results for other quantile levels are presented in Section S.3.2, Figures S2 and S3 of the Supplementary material.



The experimental results in Figure~\ref{fig:plot_compare} demonstrate that our GBD algorithm consistently outperforms CNLS-A in both computational efficiency and solution quality, offering a scalable solution to the challenges raised in \citet{dai2023variable}. By avoiding the repeated solution of large-scale integer programs with partial convexity constraints, GBD achieves significant computational savings. These benefits are particularly evident in high-noise settings.

\subsubsection{Comparison with the sparse convex regression}
As noted in Section~\ref{sec:Intro}, quantile regression extends mean regression by modeling the conditional distribution at different quantile levels. Here, we construct a dataset where the relevant features vary across quantiles to test whether our algorithm can identify the true features for each level. In contrast, sparse convex regression, such as the dual cutting-plane method in \citet{bertsimas2021sparse} that estimates only the conditional mean, is expected to fail, as it can recover only features relevant to the mean.

We generate the feature matrix $\mathbf{X}$ as described in Section~6.2.1. The true support set $S^*$, of size 10, is given by $\{0, 1, 4, 7, 8, 12, 14, 18, 24, 25\}$, where $d$ is set to either $30$ or $50$ to represent different scenarios. Motivated by \cite{lee2014model}, the response variable $y_i$ is generated according to the following function:
\[
y_i = 
\underbrace{x_0^2 + x_1^2 + x_2^2 + x_7^2}_{S^*_{\text{median}}} + 
\underbrace{\left(x_8^2 + x_{12}^2 + x_{14}^2 + x_{18}^2 + x_{24}^2 + x_{25}^2\right)}_{S^* \setminus S^*_{\text{median}}}
\cdot \epsilon \tag{17}
\]
where $\epsilon \sim \mathcal{N}(0, \sigma)$ denotes the normally distributed noise term. A notable property of this data-generating process is that the true support set for median regression (i.e., at quantile level $0.5$) is  $ S_{median}^* = \{0, 1, 4, 7\}$, while for quantile levels above $0.5$, the true active variables correspond to the full set $S^*$. We will also examine the false discovery rate (FDR) of the estimator as a complementary measure to accuracy. We tune the key parameters \(k\), \(C\), and \(\varepsilon\) via five-fold cross-validation. Specifically, we consider \(k \in \{4, \ldots, 12\}\), \(C \in \{0.1, 1, 10, 100\}\), and \(\varepsilon \in \{0, 0.02, 0.2\}\), following established practices in the literature \citep{bertsimas2021sparse, dai2023variable}.

The comparative results between our GBD algorithm and the dual cutting-plane method proposed in \cite{bertsimas2021sparse} are reported in Tables~\ref{tab:CQRvsCR_acc} and~\ref{tab:CQRvsCR_fdr}. To ensure model validity, we restrict the experiments to quantile levels \(\tau = 0.5\) and \(\tau = 0.75\), under which the conditional quantile functions remain convex. For \(\tau < 0.5\), the convexity assumption is generally violated, rendering the SCQR model inapplicable.

\begin{table}[ht]
\centering
\caption{Accuracy ($\%$) comparison of our GBD algorithm and dual cutting-plane (DCP) method in \cite{bertsimas2021sparse} (n = 1000)} 
\label{tab:CQRvsCR_acc}
\begin{threeparttable}
\resizebox{0.9\linewidth}{!}{
\begin{tabular}{ccccccccccc}
\toprule
$k$ & $d$ & $\rho$ 
& \multicolumn{4}{c|}{GBD} 
& \multicolumn{4}{c}{DCP} \\
& & 
& \multicolumn{2}{c|}{$\sigma = 0.5$} 
& \multicolumn{2}{c|}{$\sigma = 1$}
& \multicolumn{2}{c|}{$\sigma = 0.5$}
& \multicolumn{2}{c}{$\sigma = 1$} \\
& & 
& $\tau = 0.5$ & $\tau = 0.75$ 
& $\tau = 0.5$ & $\tau = 0.75$ 
& $\tau = 0.5$ & $\tau = 0.75$ 
& $\tau = 0.5$ & $\tau = 0.75$ \\
\midrule
10 & 30 & 0.2 & 95 & \textbf{94} & 95 & \textbf{94} & 100 & 40 & 95 & 40 \\
   &    & 0.5 & 90 & \textbf{90} & 85 & \textbf{88} & 95 & 38 & 90 & 38 \\
   & 50 & 0.2 & 90 & \textbf{88} & 90 & \textbf{90} & 95 & 40 & 90 & 38 \\
   &    & 0.5 & 85 & \textbf{82} & 80 & \textbf{80}  &85& 36 & 80 &  36\\
\bottomrule
\end{tabular}
}
\begin{tablenotes}
\footnotesize
\item Note: At $\tau=0.5$, accuracy is computed against $S^*_{\text{median}}$; at $\tau=0.75$, it is computed against $S^*$.
\end{tablenotes}
\end{threeparttable}
\end{table}

\begin{table}[ht]
\centering
\caption{FDR ($\%$) comparison of our GBD algorithm and DCP method in \cite{bertsimas2021sparse} (n = 1000)} 
\label{tab:CQRvsCR_fdr}
\resizebox{0.9\linewidth}{!}{
\begin{tabular}{ccccccccccc}
\toprule
$k$ & $d$ & $\rho$ 
& \multicolumn{4}{c|}{GBD} 
& \multicolumn{4}{c}{DCP} \\
& & 
& \multicolumn{2}{c|}{$\sigma = 0.5$} 
& \multicolumn{2}{c|}{$\sigma = 1$}
& \multicolumn{2}{c|}{$\sigma = 0.5$}
& \multicolumn{2}{c}{$\sigma = 1$} \\
& & 
& $\tau = 0.5$ & $\tau = 0.75$ 
& $\tau = 0.5$ & $\tau = 0.75$ 
& $\tau = 0.5$ & $\tau = 0.75$ 
& $\tau = 0.5$ & $\tau = 0.75$ \\
\midrule
10 & 30 & 0.2 & 0 & 0 &  0& 0 & 0 & 0 & 0 & 0 \\
   &    & 0.5 & 0 & 0 & 0 & 2 & 0 & 0 & 0 & 0 \\
   & 50 & 0.2 & 0 & 0 & 0 & 0 & 0 & 0 & 0 &0  \\
   &    & 0.5 & 0 & 2 & 0 & 2 & 0 &0  & 0 & 0 \\
\bottomrule
\end{tabular}
}
\end{table}

We can see that GBD algorithm can identify most true variables at different quantile levels (with an accuracy over $80\%$), while the DCP algorithm can only identify the true variables at mean value (at $\tau = 0.75$, the accuracy is less than $40\%$.). These results underscore the advantages of the SCQR framework and highlight that our GBD algorithm provides a promising approach to addressing the computational challenges raised in \cite{dai2023variable}.

\subsection{Experiments with real data}
To demonstrate the practical value of our proposed scalable algorithm, we apply it to the real-world Sustainable Development Goals benchmarking problem introduced in \cite{dai2023variable}, where the SCQR model was originally proposed. The dataset, sourced from \cite{sachs2017sdg, sachs2022sustainable}, includes 25 SDG indicators for 35 OECD countries, with various social, economic, and environmental factors as inputs, and GDP growth as the output. Following the same setup as in \cite{dai2023variable}, we estimate the quantile production function using panel data from 2017, 2019, and 2020, yielding a total of 105 observations. A complete list of input variables and their descriptions is provided in Section S.3.2, Table S2 of the Supplementary material.

The SCQR model is particularly well-suited for this task, as it allows for the estimation of conditional quantiles of GDP growth based on multidimensional SDG inputs. This capability enables the construction of a series of development frontiers corresponding to different quantile levels, offering a nuanced view of country performance. This raises a natural question: why not simply use GDP to rank countries? While GDP provides a useful measure of economic output, it fails to capture the sustainability or efficiency of development. A country may achieve high GDP growth at the expense of severe environmental degradation—such as excessive $SO_2$ emissions. Furthermore, the SCQR approach offers actionable policy insights. By identifying which SDG-related factors most influence a country’s position relative to the development frontier, the model can guide targeted interventions. For further details, we refer the reader to \cite{dai2023variable}, which provides a comprehensive analysis of using SCQR to benchmark the degree of SDG achievement among OECD countries, and illustrates how the results can be utilized to guide policy implementation and inform resource allocation strategies. Therefore, the primary goal of this example is to demonstrate the practical applicability of our algorithm to real-world problems. 

For the SDG application, we also use the 5-fold cross validation procedure to determine the optimal tuning parameters $k$, $C$, and $\epsilon$. Specifically, $k$ is selected from the range $[1, 24]$, $C$ is chosen from\( \{0.1, 1, 10, 100\} \), and $\varepsilon$ is selected from $\{0, 0.4, 0.8, 1\}$. As in \cite{dai2023variable}, we evaluate our algorithm at quantile levels $\tau = 0.05, 0.35, 0.65, 0.95$, where SCQR identifies different true variables based on the cross-validated quantile loss.


\begin{table}[ht]
\centering
\caption{Performance comparison between GBD and CNLS-A algorithms}
\label{tab:perform_compare}

\scriptsize
\setlength{\tabcolsep}{3pt}
\renewcommand{\arraystretch}{1.12}

\begin{subtable}{0.49\linewidth}
\centering
\caption{Performance by GBD algorithm}
\scalebox{\tablescale}{%
\begin{tabular}{
  >{\centering\arraybackslash\rowstrut}m{0.85cm}
  >{\centering\arraybackslash\rowstrut}m{1.10cm}
  >{\centering\arraybackslash\rowstrut}m{1.25cm}
  >{\centering\arraybackslash\rowstrut}m{0.55cm}
  >{\centering\arraybackslash\rowstrut}m{0.90cm}
  >{\centering\arraybackslash\rowstrut}m{1.15cm}
  >{\centering\arraybackslash\rowstrut}m{1.15cm}
}
\toprule
$\tau$ & MAE (in) & MAE (out) & $k$ & time (s) & F-MAE (in) & F-MAE (out) \\
\midrule
0.05 & 0.79 & \textbf{2.39} & 8 & 287 & 0.50 & 2.65 \\
0.35 & 3.26 & \textbf{3.66} & 4 & 389 & 1.59 & 4.83 \\
0.65 & 3.38 & \textbf{4.38} & 3 & 400 & 1.89 & 4.49 \\
0.95 & 1.04 & \textbf{1.15} & 6 & 98  & 0.72 & 1.36 \\
\bottomrule
\end{tabular}}
\end{subtable}%
\hfill
\begin{subtable}{0.49\linewidth}
\centering
\caption{Performance by CNLS-A algorithm}
\scalebox{\tablescale}{%
\begin{tabular}{
  >{\centering\arraybackslash\rowstrut}m{0.85cm}
  >{\centering\arraybackslash\rowstrut}m{1.35cm}
  >{\centering\arraybackslash\rowstrut}m{1.35cm}
}
\toprule
$\tau$ & MAE (in) & MAE (out) \\
\midrule
0.05 & 0.38 & 2.60 \\
0.35 & 3.09 & 4.04 \\
0.65 & 3.25 & 4.41 \\
0.95 & 0.96 & 1.22 \\
\bottomrule
\end{tabular}}
\end{subtable}

\end{table}

\noindent
Table \ref{tab:perform_compare} summarizes the performance comparison results of the GBD and CNLS-A algorithms across various quantile levels. \textbf{MAE (in/out)} represent the quantile losses on the training/testing sets, and \textbf{$k$} is the number of selected relevant variables. For comparison, \textbf{F-MAE (in/out)} report the losses without variable selection. These results demonstrate that our GBD algorithm achieves better generalization and more compact models than CNLS-A and the no-selection baseline. Detailed regression results (selected features along with their average estimated coefficients) in Section S.3.2, Table S.3 of the Supplementary material further show distinct model structures across quantile levels.

\section{Discussion}
Subset selection in high-dimensional settings remains a challenging NP-hard problem. This paper proposes a scalable GBD framework for SCQR, where the subproblem is efficiently solved via an adapted cutting-plane method. To accelerate convergence, we further incorporate a warm-start strategy and a novel local search-based matheuristic.

Extensive experiments validate the effectiveness of our framework. Compared to standard CQR and Lipschitz-constrained models, our regularized formulation offers improved generalization. Against the CNLS-A algorithm \citep{dai2023variable}, GBD achieves notable gains in both runtime and variable selection accuracy. A real-world application to SDG evaluation across OECD countries further demonstrates its practical value. Beyond computational performance, our SDG case study shows that SCQR uncovers heterogeneity in development achievements across OECD countries. Such benchmarking enables cross-country policy comparison and supports evidence-based prioritization of indicators, offering practical guidance for tailored interventions and resource allocation.

Future work will focus on strengthening cut generation and extending the LSB framework, as well as analyzing the convergence rate of the decomposition algorithm under certain structural conditions. Another avenue is to explore limiting the number of pieces in SCQR models, aiming to reduce computational complexity while preserving the accuracy of variable selection.

\section*{Acknowledgments}
This work was funded by the National Nature Science Foundation of China under Grant No. 12320101001 and 12071428.

\section*{Disclosure statement}
Declarations of interest: none





\bibliographystyle{model5-names} 
\bibliography{ref.bib}






\end{document}